\documentclass[11pt]{article}%
\usepackage{amssymb}
\usepackage{amsmath}
\usepackage{amsfonts}
\usepackage{graphicx}%
\newtheorem{theorem}{Theorem}[section]

\newtheorem{proposition}{Proposition}[section]
\newtheorem{remark}{Remark}[section]

\newenvironment{proof}[1][Proof]{\noindent\textbf{#1.} }{\ \hfill\rule{0.3em}{0.5em}}

\newcommand*{\C}{\mathbb{C}}

\def\beq#1{\begin{equation}\label{#1}}
\def\eeq{\end{equation}}
\def\bep{\begin{proof}}

\def\ep{\end{proof}}

\def\y{\mathbf y}

\def\A{{\mathcal A}}
\def\B{\mathcal B}
\def\ST{{\mathcal T}}
\def\OO{{\mathcal O}}
\def\ignore#1{}

\def\qed{\hfill {$ \Box $} \medskip}
\def\bea#1{\begin{array}{#1}}
\def\ea{\end{array}}

\input epsf

\begin{document}

\title{The Quantum Eigenvalue Problem and Z-Eigenvalues of Tensors}
\author{Xinzhen Zhang\footnote{\small Department of
 Mathematics, School of Science, Tianjin University, Tianjin,
 300072, China. Her work is supported by the National Natural Science Foundation of China (Grant No. 11101303
and 11171180)
 E-mail: \emph{xzzhang@tju.edu.cn}.}\quad  and\quad
Liqun Qi\thanks{\small Department of Applied Mathematics, The Hong
Kong Polytechnic University, Hung Hom, Kowloon, Hong Kong.  E-mail:
\emph{maqilq@polyu.edu.hk}. His work is supported by the Hong Kong
Research Grant Council. }}

\date{\today}
\maketitle
\begin{quote}
{\small \noindent\textbf {Abstract.} The quantum eigenvalue problem
arises in the study of the geometric measure of the quantum
entanglement.   In this paper, we convert the quantum eigenvalue
problem to the Z-eigenvalue problem of a real symmetric tensor.  In
this way, the theory and algorithms for Z-eigenvalues can be applied
to the quantum eigenvalue problem.  In particular, this gives an
upper bound for the number of quantum eigenvalues.    We show that
the quantum eigenvalues appear in pairs, i.e., if a real number
$\lambda$ is a quantum eigenvalue of a square symmetric tensor
$\Psi$, then $-\lambda$ is also a quantum eigenvalue of $\Psi$. When
$\Psi$ is real, we show that the entanglement eigenvalue of $\Psi$
is always greater than or equal to the Z-spectral radius of $\Psi$,
and that in several cases the equality holds. We also show that the
ratio between the entanglement eigenvalue and the Z-spectral radius
of a real symmetric tensor is bounded above in a real symmetric
tensor space of fixed order and dimension. }
\\
{\small \medskip\noindent\textbf {Key Words.}  quantum eigenvalue,
Z-eigenvalue, real symmetric tensor, entanglement eigenvalue.}
\end{quote}
\section{Introduction}\label{Sec1}
Eigenvalues of higher order tensors were introduced in 2005
\cite{qi, li} and have attracted much attention in the literature.
The E-eigenvalues and E-eigenvectors of tensors were introduced in
\cite{qi, qi07}.    For a real tensor, if an E-eigenvalue has a real
E-eigenvector, then it is also real and called a Z-eigenvalue of
that tensor.    The corresponding E-eigenvector is called a
Z-eigenvector. Z-eigenvalues found applications in determining
positive definiteness of a multivariate system \cite{qi, qww}, best
rank-one approximation \cite{qi, qi11, qww, zqy}, magnetic resonance
imaging \cite{bv, qyw}, spectral hypergraph theory \cite{hq, lqy}
and symmetric embedding \cite{rv}, and possess links with geometry
\cite{ba, ba09, ba11, ba10, qi06}.     Properties of and algorithms
for finding Z-eigenvalues and Z-eigenvectors can be found in
\cite{cpz, cpz1, cz, km, qww, qyw, za}. Z-eigenvalues are roots of
E-characteristic polynomials, which have been studied in \cite{cs,
lqz, nqww, qi, qi07}.

Theoretically, people were not very sure if Z-eigenvalues or E-eigenvalues are more essential.   Actually, the definition of E-eigenvalues missed some eigenvalues in the matrix case \cite{km, qi07}.   Another concept, equivalence
eigenpairs were introduced to remedy this defect in \cite{cs}.   People were also not sure if the normalizing constraint in the definition of Z-eigenvalues has any application background.

The quantum eigenvalue problem arises in the study of the
geometric measure of the quantum entanglement \cite{wg}.   The
quantum entanglement problem is a central problem in quantum
information \cite{nc}.  The geometric measure of entanglement of a
multipartite pure state, based upon the minimum distance of that
state from the set of separable pure states, was introduced in
\cite{wg} and has been studied in
\cite{bh,cxz,hkwg,hs,hmmov,odv,qi12}.   The quantum eigenvalues are
real numbers.    The largest quantum eigenvalue is called the {\bf
entanglement eigenvalue}, which is corresponding to the nearest separable
pure state to the given pure state, and thus has its physical meaning.

Only recently, it was shown in \cite{hqz} that the entanglement
eigenvalue of a symmetric pure state with nonnegative amplitudes is
equal to the largest Z-eigenvalue of the underlying nonnegative
tensor.

In this paper, we convert the quantum eigenvalue problem to a
Z-eigenvalue problem of a real symmetric tensor.  In this way, on
one hand, the theory and algorithms for Z-eigenvalues can be applied
to the quantum eigenvalue problem. For example, we now know that the
number of quantum eigenvalues of a multipartite pure state is
finite, and we know an upper bound of that number.  These were
previously unknown. On the other hand, we now know Z-eigenvalues are
more essential, E-eigenvalues and equivalence eigenpairs are some
useful concepts for studying Z-eigenvalues, and the normalizing
constraint for Z-eigenvalues has some physics background in this
case.

We show that the quantum eigenvalues appear in pairs, i.e., if a real number $\lambda$ is a quantum eigenvalue of a square symmetric tensor $\Psi$, then $-\lambda$ is also a quantum eigenvalue of $\Psi$.    When $\Psi$ is real, we show that the entanglement eigenvalue of $\Psi$ is always greater than or equal to the Z-spectral radius of $\Psi$, and that in several cases the equality holds.

The rest of the paper is organized in the following way.   In the
next section, the Z-eigenvalue problem, the quantum eigenvalue
problem and some related concepts are reviewed.  We show there that
the quantum eigenvalues appear in pairs.      In Section 3, we show
that when the underlying symmetric tensor $\Psi$ is real, the
entanglement eigenvalue of $\Psi$ is always greater than or equal to
the Z-spectral radius of $\Psi$, and that in several cases the
equality holds. In general, the equality may not hold.  Then, for
the $m$-order $n$-dimensional real symmetric tensor space, we
introduce the QR-ratio which is the maximum ratio between the
entanglement eigenvalues and the Z-spectral radii of the tensors in
that space.    We convert the quantum eigenvalue problem to the
Z-eigenvalue problem of a real symmetric tensor, and discuss the
consequences of this result in Section 4. We also show in Section 4
that the QR-ratio is finite there. Some final remarks are made in
Section 5.

Throughout this paper, we assume that $m, n \ge 2$.

\section{The Z-eigenvalue problem and the quantum eigenvalue
problem}

Suppose that $\ST = (\ST_{i_1\cdots i_m})$ is an $m$-order
$n$-dimensional tensor, where $\ST_{i_1\cdots i_m} \in \C$ for $i_1,
\cdots, i_m = 1, \cdots, n$. If there are a number $\lambda \in \C$
and a nonzero vector $w \in \C^n$ such that
\begin{equation}\label{Zeig} \left\{\begin{array}{rl}\ST w^{m-1}&=\lambda w,\\
w^\top w&=1, \end{array}\right.
\end{equation}
where $\ST w^{m-1} \in \C^n$ with its $i$th component defined by
$$\left(\ST w^{m-1}\right)_i = \sum\limits_{i_2, \cdots, i_m=1}^m
\ST_{ii_2\cdots i_m}w_{i_2}\cdots w_{i_m},$$ then we call $\lambda$
an {\bf E-eigenvalue} of $\ST$ with $w$ as its corresponding {\bf
E-eigenvector}.   We always have
$$\lambda = \ST w^m \equiv \sum_{i_1, \cdots, i_m=1}^m
\ST_{i_1\cdots i_m}w_{i_1}\cdots w_{i_m}.$$ If $\ST$ and $w$ are
real, then $\lambda$ is also real.  In this case, we call $\lambda$
a {\bf Z-eigenvalue} of $\ST$ with $w$ as its corresponding {\bf
Z-eigenvector}.  By \cite{qi}, if $\ST$ is real and symmetric,
Z-eigenvalues always exist.  In \cite{qi11}, the largest absolute
value of Z-eigenvalues of $\ST$ is called the spectral radius of
$\ST$. In \cite{cpz1}, this value is called the {\bf Z-spectral
radius} of $\ST$.   Since there are other kinds of eigenvalues of
tensors, we use the name Z-spectral radius in this paper.  We denote
the Z-spectral radius of $\ST$ as $Z(\ST)$.       For more
properties of Z-eigenvalues, we refer to the references cited in the
first paragraph of this paper.

Define the Frobenius norm of $\ST$ as
$$\| \ST \| = \sqrt{ \sum_{i_1, \cdots, i_m = 1}^n \ST_{i_1\cdots i_m}^2}.$$

We say that tensor $\ST$ is {\bf symmetric} if its entries
$\ST_{i_1\cdots i_m}$ are invariant under any permutation of its
indices.

Denote the space of $m$-order $n$-dimensional real
symmetric tensors by $S(m, n)$.  Denote the zero tensor in $S(m, n)$ by $\OO$.  It was proved in \cite{qi11} that the Z-spectral radius is also a norm of $S(m, n)$.  By the norm equivalence theorem in finite dimensional spaces \cite{or}, both
$$\rho(m, n) := \sup \left\{ {Z(\ST) \over \|\ST\|} : \ST \in S(m, n) \setminus \{ \OO \} \right\}$$
and
$$\mu(m, n) := \inf \left\{ {Z(\ST) \over \|\ST\|} : \ST \in S(m, n) \setminus \{ \OO \} \right\}$$
are finite positive numbers.  In \cite{qi11}, $\mu(m, n)$ is called the {\bf best rank-one approximation ratio} of $S(m, n)$.

In this paper, we only discuss the quantum eigenvalue problem
arising from the geometric measure of entanglement of a multipartite
symmetric pure state.   The quantum eigenvalue problem for the
geometric measure of entanglement of a multipartite nonsymmetric
pure state can be converted to the symmetric one by symmetric
embedding \cite{rv}, see \cite{hqz}.

 An $m$-partite symmetric pure state $| \Psi \rangle$ can be regarded as an element in a tensor product space ${\mathcal H}^{\otimes m}$, satisfying $\langle \Psi | \Psi \rangle = 1$, where the dimension of $\mathcal H$ is $n$.  An $m$-partite symmetric pure separable state has the form $|\phi \rangle^{\otimes m}$, where $|\phi \rangle \in {\mathcal H}$ and $\langle \phi | \phi \rangle = 1$.    Then, by \cite{wg, hkwg, qi12}, the quantum eigenvalue problem has the form
 \begin{equation}\label{Qeig} \left\{\begin{array}{rl}\langle \Psi | \phi \rangle^{\otimes (m-1)}&=\lambda \langle \phi |,\\
\langle \phi^{\otimes (m-1)} | \Psi \rangle &=\lambda | \phi \rangle,\\
\langle \phi | \phi \rangle &=1, \\
\langle \Psi | \Psi \rangle &=1.\end{array}\right.
\end{equation}

 We assume that  under an orthonormal basis of $\mathcal H$, by abusing the notation, the state $| \Psi \rangle$ is represented by an $m$-order $n$-dimensional symmetric complex tensor $\Psi = (\Psi_{i_1\cdots i_m})$,  and the state $| \phi \rangle$ is represented by a
vector $z \in \C^n$.    The assumption that $\langle \phi | \phi \rangle = 1$ implies that $\bar z^\top z = 1$, where the bar denotes the conjugate.   We keep this assumption.   The assumption that $\langle \Psi | \Psi \rangle = 1$ implies that the Frobenius norm of $\Psi$ should be one. We
drop this assumption to accommodate more general problems.   Throughout this paper, we assume that $\Psi$ is symmetric.

With above notation, the quantum eigenvalue problem (\ref{Qeig}) has
the following form.
\begin{equation}\label{eig} \left\{\begin{array}{rl}\Psi z^{m-1}&=\lambda \bar z,\\
\bar\Psi \bar z^{m-1}&=\lambda z,\\
\bar z^\top z&=1.\end{array}\right.
\end{equation}
Suppose $\lambda \in \C$ and $z \in \C^n$ satisfy (\ref{eig}). Then
$\lambda$ must be real.  We call it a {\bf quantum eigenvalue} or in
short a {\bf Q-eigenvalue} of the tensor $\Psi$, and call $z$ a {\bf
Q-eigenvector} of $\Psi$, associated with the Q-eigenvalue
$\lambda$.   The largest Q-eigenvalue of $\Psi$ is called the {\bf
entanglement eigenvalue} of $\Psi$, and denoted as $Q(\Psi)$.   It
is associated with the nearest separable pure state of the state
$|\Psi \rangle$.

In fact, (\ref{eig}) is equivalent to
\begin{equation}\label{eig1}  \left\{\begin{array}{rl}\Psi z^{m-1}&=\lambda \bar z,\\
\bar z^\top z&=1, \\~~~~\lambda&\in \Re.\end{array}\right.
\end{equation}

By \cite{qi12, wg}, Q-eigenvalues always exist.

Because of the conjugate operation is involved in (\ref{eig}) and
(\ref{eig1}), the theory of algebraic geometry \cite{clo, gkz}
cannot be applied to them directly.  This poses difficulty to
analyze them. Before this paper, it was unknown even if there are a
finite number of Q-eigenvalues or not.

We now prove a property of Q-eigenvalues.

\begin{proposition} \label{p1}
The Q-eigenvalues appear in pairs, i.e., if a real number $\lambda$ is a  Q-eigenvalue of a square symmetric tensor $\Psi$, then $-\lambda$ is also a Q-eigenvalue of $\Psi$.
\end{proposition}
\proof Suppose that $\lambda$ is a Q-eigenvalue of $\Psi$ with its corresponding Q-eigenvector $z$.   Then $-\lambda$ is a Q-eigenvalue of $\Psi$ with its corresponding Q-eigenvector $ze^{\pi\sqrt{-1} \over m}$.  The conclusion follows.  \qed

\section{The case that $\Psi$ is real}

Recently, it was proven in \cite{hqz} that if furthermore $\Psi$ is
nonnegative, then $Q(\Psi) = Z(\Psi)$.    We now establish results
in more general cases.

\begin{theorem}\label{rthm} If $\Psi$ is real, then $Q(\Psi) \ge Z(\Psi)$.   In the following six cases, we have $Q(\Psi) =
Z(\Psi)$.   We assume that $\Psi$ is real in all these six cases.

1). $m = 2$;

2). $\Psi$ is diagonal;

3). $\Psi$ is nonnegative;

4). $\Psi$ is nonpositive;

5). $\Psi=\sum\limits_{k=1}^n\alpha_k \left({y^{(k)}}\right)^m$
where $\alpha_k, k= 1, \cdots n$, are real numbers, $\{ y^{(k)} : k
= 1, \cdots, n \}$ is an orthonormal basis of $\Re^n$;

6). $m \ge 4$ is even, an element of $\Psi$ is nonzero only if a
half of its indices are the same and the other half of its indices
are also the same, and $\Psi$ is diagonally dominated in the sense
that the absolute value of each of its diagonal elements is greater
than or equal to $\left({m-1 \atop {m \over 2}}\right)$ times of the
absolute value of any off-diagonal element whose indices overlap
with that diagonal element.
\end{theorem}
\proof By the definitions of Q-eigenvalues and Z-eigenvalues, we see
that if $\Psi$ is real, each Z-eigenvalue of $\Psi$ is also a
Q-eigenvalue of $\Psi$. By Proposition \ref{p1}, if $\lambda$ is a
Q-eigenvalue of $\Psi$, then $-\lambda$ is also a Q-eigenvalue of
$\Psi$.  Putting these together, we have $Q(\Psi) \ge Z(\Psi)$ when
$\Psi$ is real.

1). This case can be regarded as a special case of case 5). Here we
give a direct proof for this case.  Assume that $\lambda$ is the
entanglement eigenvalue of $\Psi$ with $z$ as its corresponding
Q-eigenvector. Assume that $z= x + y\sqrt{-1}$, where $x$ and $y$
are real.   Then we have $\Psi x = \lambda x$, $\Psi y = - \lambda
y$ and at least one of $x$ and $y$ is not equal to zero. If $x \not
= 0$, then $\lambda$ is a Z-eigenvalue of $\Psi$ with $x \over
\sqrt{x^\top x}$ as its corresponding Z-eigenvector. If $y \not =
0$, then $-\lambda$ is a Z-eigenvalue of $\Psi$ with $y \over
\sqrt{y^\top y}$ as its corresponding Z-eigenvector.   This shows
that $Q(\Psi) \le Z(\Psi)$. Since we always have $Q(\Psi) \ge
Z(\Psi)$, we have $Q(\Psi) = Z(\Psi)$.

2). Suppose that $\Psi$ is real and diagonal.  Denote the diagonal
elements of $\Psi$ by $a_{k\cdots k}$ for $k= 1, \cdots, n$. Assume
that $|a_{i\dots i}| = \max \{ |a_{k\cdots k}| : k = 1, \cdots, n
\}$.   Let $x_i = 1$ and $x_k = 0$ for $k \not = i$. Let $\mu =
a_{i\cdots i}$.  Then $\mu$ is a Z-eigenvalue of $\Psi$ with $x$ as
its corresponding Z-eigenvector.   This shows that $Z(\Psi) \ge
|a_{i\cdots i}|$.  Let $\lambda$ be the entanglement eigenvalue of
$\Psi$ with $z$ as its corresponding Q-eigenvector. Assume that $z_j
\not = 0$.   Then $|z_j| \le \sqrt{\bar z^\top z} = 1$, and we have
$a_{j\cdots j}z_j^{m-1} = \lambda \bar z_j$, which implies that
$$Q(\Psi) = \lambda \le |a_{j\cdots j}|\cdot |z_j|^{m-2} \le |a_{i\cdots i}| \le Z(\Psi).$$    Since we always have $Q(\Psi) \ge Z(\Psi)$, we have $Q(\Psi) = Z(\Psi)$.

3). This was proved in \cite{hqz}.

4). By definition, $Z(\Psi) = Z(-\Psi)$.  By Proposition \ref{p1}, $Q(\Psi) = Q(-\Psi)$.   By 3), $Q(-\Psi) = Z(-\Psi)$.  The conclusion follows now.

5). It suffices to prove that in this case, if $\lambda$ is a
Q-eigenvalue of $\Psi$, then either $\lambda$ or $-\lambda$ is a
Z-eigenvalue of $\Psi$.    Suppose  that $\lambda$ is a Q-eigenvalue
of $\Psi$, with Q-eigenvector $z$. Then $\Psi z^{m-1}=\lambda \bar
z$. That is,
$$\sum\limits_{i=1}^n\alpha_k ({y^{(k)}}^\top z)^{m-1} y^{(k)}=\lambda \bar z.$$

For convenience of notation, we denote $\beta_k={y^{(k)}}^\top z=r_k
e^{\theta_k\sqrt{-1}}$ for $k=1,2,\cdots, n$. Since $y^{(1)},
\cdots, y^{(n)}$ are orthonormal vectors, there exists $k$ such that
$\beta_k={y^{(k)}}^\top z=r_k e^{\theta_k\sqrt{-1}}\neq 0$.

Multiplying $(y^{(k)})^\top$ on the both sides of above equality, we
have $\alpha_k r_k^{m-2}e^{m\theta_k \sqrt{-1}}=\lambda$, which
implies that $m\theta_k=l\pi$ for an integer $l\geq 0$.   Hence
$\lambda=\alpha_k r_k^{m-2}\cos{l\pi} $ for all $k$ satisfying
$\beta_k\neq 0$.

Let $x=\sum\limits_{j=1}^nr_j y^{(j)}\in\Re^n$. Then $x^\top x=1$
and
$$\Psi x^{m-1}=\sum\limits_{k=1}^n \alpha_ky^{(k)}r_k^{m-1}=\sum\limits_{k=1}^n y^{(k)}(-1)^l\lambda r_k=(-1)^l\lambda x,$$
which means that $(-1)^l \lambda$ is a Z-eigenvalue of $\Psi$.  The
proof of 5) is completed.

6). Suppose that $\Psi$ satisfies the conditions. Denote the
elements of $\Psi$ by $a_{k_1\cdots k_m}$ for $k_1, \cdots , k_m =
1, \cdots, n$.    Assume that $|a_{i\cdots i}| = \max \{ |a_{k\cdots
k}| : k = 1, \cdots, n \}$.   Let $x_i = 1$ and $x_k = 0$ for $k
\not = i$. Let $\mu = a_{i\cdots i}$.  Then $\mu$ is a Z-eigenvalue
of $\Psi$ with $x$ as its corresponding Z-eigenvector, as an element
of $\Psi$ is nonzero only if a half of its indices are the same and
the other half of its indices are also the same.   This shows that
$Z(\Psi) \ge |a_{i\cdots i}|$.  Let $\lambda$ be the entanglement
eigenvalue of $\Psi$ with $z$ as its corresponding Q-eigenvector.
Assume that $z_j \not = 0$.   Then $|z_j| \le \sqrt{\bar z^\top z} =
1$, and we have $a_{j\cdots j}z_j^{m-1} + \sum_{k \not = j}
\left({m-1 \atop {m \over 2}}\right)a_{j\cdots jk\cdots k}z_j^{m-2
\over 2}z_k^{m \over 2} = \lambda \bar z_j$, which implies that
$$Q(\Psi) = \lambda \le |a_{j\cdots j}|\cdot \sum_{k=1}^n |z_k|^{m \over 2} \le |a_{i\cdots i}| \le Z(\Psi).$$    Since we always have $Q(\Psi) \ge Z(\Psi)$, we have $Q(\Psi) = Z(\Psi)$.
\qed

In general, the equality $Q(\Psi) \ge Z(\Psi)$ may not hold. On Page
3 of \cite{hkwg}, a counter example is given.   By computation, for
that example, we have $Z(\Psi) = 2.2805/\sqrt{21}$ while $Q(\Psi) =
3.1768/\sqrt{21}$.

If the equality $Q(\Psi) \ge Z(\Psi)$ holds, then the problem for
finding the entanglement eigenvalue of a real symmetric tensor
$\Psi$ is turned out to the problem for finding the Z-spectral
radius of $\Psi$.   This reduces a half of the number of variables.
Thus, a challenging task is to find a sufficient and necessary
condition or some more general sufficient conditions such that this
equality holds.

Another challenging task is to find the maximum ratio of $Q(\Psi)$
and $Z(\Psi)$.      Define the {\bf QR-ratio} of $S(m,
n)$ as
$$R(m, n) = \sup \left\{ { Q(\Psi) \over Z(\Psi)} : \Psi \in S(m, n)
\setminus \{ \OO \} \right\}.$$ By Theorem \ref{rthm}, we have $R(2,
n) = 1$ and $R(m, n) \ge 1$ for $m \ge 3$.  In the next section, we
will show that the QR-ratio $R(m, n)$ is a finite positive number.
Even if we know that it is a finite positive number, what is its
value, and is it attainable?  These remain as topics for further
research.

\section{Converting the quantum eigenvalue problem to the Z-eigenvalue problem of a real symmetric tensor}\label{se2}

In this section, we convert the quantum eigenvalue problem to the Z-eigenvalue problem of a real symmetric tensor.

Before proceeding, the following relation is established when $m=2$.
\begin{theorem}\label{ref1} Suppose that $m=2$. Then there exists a symmetric matrix $M\in \Re^{2n\times 2n}$ such that a number $\lambda$ is a
Q-eigenvalue of the tensor $\Psi$ if and only if it is a Z-eigenvalue
of the real symmetric matrix $M$.   In this case, $z = x +
y\sqrt{-1}$ is a Q-eigenvector of $\Psi$, if and only if $w=({x
\atop y})$ is a unit eigenvector of $M$, associated with $\lambda$.

Furthermore, the eigenvalues of $M$ appear in pairs, i.e., if $\lambda$ is an eigenvalue of $M$, then $-\lambda$ is also an eigenvalue of $M$.
\end{theorem}

\proof For convenience of notation, let $\Psi=A+B\sqrt{-1}$ and
$z=x+y\sqrt{-1}$ be a Q-eigenvector of $\Psi$ associated with the
Q-eigenvalue $\lambda$. Then $\Psi z=\lambda \bar z$ and $\bar
z^\top z=1$, that is
$$\left\{\begin{array}{rl}Ax-By=\lambda x,\\-Bx-Ay=\lambda y,\\ w^\top w =1.\end{array}\right.$$

Let $M=\left(\begin{array}{rl}A& -B\\-B &-A\end{array}\right)$ and
$w=\left(\begin{array}{rl}x\\y\end{array}\right)$.  Then $M$ is a
real symmetric matrix, and the above system holds if and only if
$$\left\{\begin{array}{rl}M w=\lambda w, \\
w^\top w = 1. \end{array}\right.$$

Furthermore, assume that $\lambda$ is an eigenvalue of $M$, with an eigenvector $w$.   Denote $w=({x \atop y})$ and $\hat w=({y \atop -x})$.   Then we see that $-\lambda$ is also an eigenvalue of $M$ with eigenvector $\hat w$.  This completes our proof.

\qed

Note that since $M$ is real symmetric, all of its eigenvalues are real, thus Q-eigenvalues of $\Psi$.

 Now we are ready to extend the result to the m-partite case with $m\geq 3$. 
 For $i = 1, \cdots, 2n$, define $\hat i = i$ if $i \le
 n$, and $\hat i = i-n$ if $i > n$.  For an $m$-order $n$-dimensional tensor
 $\ST = (\ST_{i_1\cdots i_m})$, two vectors $u = (u_i), v = (v_i) \in \C^n$ and $k = 0, \cdots, m-1$, we have $\ST u^{m-1-k}v^k \in \C^n$, defined by
 $$\left(\ST u^{m-1-k}v^k\right)_i = \sum_{i_2,\cdots, i_m = 1}^n \ST_{ii_2\cdots i_m}u_{i_2}\cdots u_{i_{m-k}}v_{i_{m-k+1}}\cdots v_{i_m}.$$

 \begin{theorem} \label{mthm} Suppose that $\Psi$ is an $m$-order
$n$-dimensional symmetric tensor, where $m\geq 3$. Then there exists
a real symmetric tensor $\ST\in \Re^{\overbrace{2n\times
\cdots\times 2n}^m}$ such that $\lambda$ is a Q-eigenvalue of $\Psi$
if and only if it is a Z-eigenvalue of $\ST$. In this case, $z = x +
y\sqrt{-1}$ is a Q-eigenvector of $\Psi$ , if and only if $w=({x
\atop y})$ is a Z-eigenvector of $\ST$, associated with $\lambda$.

Furthermore, the E-eigenvalues and Z-eigenvalues of $\ST$ appear in pairs, i.e., if $\lambda$ is an E-eigenvalue or Z-eigenvalue of $\ST$, then $-\lambda$ is also an E-eigenvalue or Z-eigenvalue of $\ST$ correspondingly.
\end{theorem}
\proof We assume that $\Psi=\A+\B \sqrt{-1}$ is an $m$-order $n$-dimensional symmetric tensor, where tensors $\A$ and $\B$ are real symmetric tensors.

Suppose that $z=x+y\sqrt{-1}$ is a Q-eigenvector of $\Psi$ associated with the
Q-eigenvalue $\lambda$.  Let $w=({x
\atop y})$. Then we have
$\Psi z^{m-1}=\lambda \bar z$ and $\bar z^\top z=1$, that is,
$$\left\{\begin{array}{rl} \left(\A+\B \sqrt{-1}\right)\sum\limits_{k=0}^{m-1}\left(\sqrt{-1}\right)^k
\left({m-1 \atop k}\right)x^{m-1-k}y^k&=\lambda \left(x -y\sqrt{-1}\right) ,\\
w^\top w&=1.\end{array}\right.$$    Considering the two cases $k=2j$
and $k=2j+1$, we have
$$\left\{\begin{array}{rl} \left(\A+\B \sqrt{-1}\right)\left[\sum\limits_{j=0}^{\left\llcorner {m-1 \over 2}\right\lrcorner}(-1)^j\left({m-1 \atop 2j}\right) x^{m-1-2j}y^{2j}+\sqrt{-1}\sum\limits_{j=0}^{\left\llcorner {m-2 \over 2}\right\lrcorner}(-1)^j\left({m-1 \atop 2j+1}\right) x^{m-2j-2}y^{2j+1}\right]&\\
=\lambda \left(x -y\sqrt{-1}\right) ,\\
w^\top w=1.&\end{array}\right.$$
Separating the real and imaginary parts, we have
$$\left\{\begin{array}{rl} \sum\limits_{j=0}^{\left\llcorner {m-1 \over 2}\right\lrcorner}(-1)^{j}\left({m-1 \atop 2j}\right)\A x^{m-1-2j}y^{2j}+\sum\limits_{j=0}^{\left\llcorner {m-2 \over 2}\right\lrcorner}(-1)^{j+1}\left({m-1 \atop 2j+1}\right)\B x^{m-2j-2}y^{2j+1} &=\lambda x,\\
\sum\limits_{j=0}^{\left\llcorner {m-2 \over 2}\right\lrcorner}(-1)^{j+1}\left({m-1 \atop 2j+1}\right)\A x^{m-2j-2}y^{2j+1}+\sum\limits_{j=0}^{\left\llcorner {m-1 \over 2}\right\lrcorner}(-1)^{j+1}\left({m-1 \atop 2j}\right)\B x^{m-1-2j}y^{2j}&=\lambda y,\\
w^\top w&=1.\end{array}\right.$$

Let tensor $\ST\in \Re^{2n\times \cdots\times 2n}$ be an $m$-order $2n$-dimensional symmetric tensor, defined by
$$\ST_{i_1i_2\cdots i_m}=\left\{\begin{array}{rl}
(-1)^j \A_{\hat i_1\hat i_2\cdots \hat i_m}&\mbox{when there exist $2j$ ($0\leq j\leq \frac{m}{2}$) indices larger than $n$},~~~~ \\
(-1)^{j+1}\B_{\hat i_1\hat i_2\cdots \hat i_m}&\mbox{when there
exist $2j+1$ ($0\leq j\leq \frac{m-1}{2}$) indices larger than
$n$}.\end{array}\right.$$ Then we have $\ST w^{m-1}=\lambda w$. This
indicates that $\lambda$ is a Q-eigenvalue of tensor $\Psi$,
associated with Q-eigenvector $z$ if and only if $\lambda$ is a
Z-eigenvalue of tensor $\ST$, associated with Z-eigenvector $w$.

Furthermore, assume that $\lambda$ is an E-eigenvalue of $\ST$, with an E-eigenvector $w=({x \atop y})$.   Here, $x$ and $y$ may not be real.  Let $\hat w=({y \atop -x})$.   We may easily see that $-\lambda$ is also an E-eigenvalue of $\ST$ with E-eigenvector $\hat w$.   Now assume that $\lambda$ is a Z-eigenvalue of $\ST$.    In this case, we may assume that $x$ and $y$ are real. Then both $w$ and $\hat w$ are real.   Thus, in this case, $-\lambda$ is also a Z-eigenvalue of $\ST$ with Z-eigenvector $\hat w$.   This completes our proof.
\qed

\begin{remark}The reformulated tensor $\ST$ is not unique. Assume that $m=3$. For example, $w=({x \atop y})$ is a Z-eigenvector of tensor  $\ST$, associated with the Z-eigenvalue $\lambda$, where $\ST$ is defined by  $$ \ST_{ijk}=\left\{\begin{array}{cl}\A_{ijk}~~& \mbox{if }~~~~ i,j,k\leq n,\\
\B_{\hat i\hat j\hat k}&\mbox{if}~~~~{\rm one \ of \ } \{ i, j, k \} >n,~~ {\rm two \ of \ } \{ i, j, k \} \le n, \\
-\A_{\hat i\hat j\hat k}&\mbox{if}~~~~{\rm one \ of \ } \{ i, j, k \} \le n,~~{\rm two \ of \ } \{ i, j, k \}>n,\\
-\B_{\hat i\hat j\hat k}&\mbox{if}~~~~i,j,k>n,\end{array}\right.$$
when $z=x+y\sqrt{-1}$ is a Q-eigenvector of $\Psi$, associated with
the Q-eigenvalue $\lambda$.
 \end{remark}

 By Theorem \ref{mthm} and the result in \cite{cs}, we now have an upper bound for the number of Q-eigenvalues.

\begin{proposition}\label{qenum}There are at most $
\frac{(m-1)^{2n}-1}{m-2}$ Q-eigenvalues for a symmetric $m$-partite pure state $|\Psi\rangle$.  \end{proposition}

 By Theorem \ref{mthm}, we now also know Z-eigenvalues are more essential, E-eigenvalues
and equivalence eigenpairs are some useful concepts for studying
Z-eigenvalues, and the normalizing constraint for Z-eigenvalues has
some physics background in this case.

We also note that $\ST$ is a special real symmetric tensor: its E-eigenvalues and Z-eigenvalues appear in pair.   Do $\ST$ and its E-characteristic polynomial have other special properties?  This leaves as a further research topic.

Finally, we prove that the QR-ratio $R(m, n)$ is finite.
\begin{theorem}
The QR-ratio $R(m, n)$ is a finite positive number.
Actually, we have
$$1 \le R(m, n) \le {2^{m-1 \over 2}\rho(m, 2n) \over \mu(m, n)}.$$
\end{theorem}
\proof The first inequality is by Theorem \ref{rthm}.   Let $\Psi
\in S(m, n)$. By Theorem \ref{mthm}, we have $Q(\Psi) = Z(\ST)$,
where $\ST$ is a real symmetric tensor in $S(m, 2n)$, defined by
$$\ST_{i_1i_2\cdots i_m}=\left\{\begin{array}{rl}
(-1)^j \Psi_{\hat i_1\hat i_2\cdots \hat i_m}&\mbox{when there exist $2j$ ($0\leq j\leq \frac{m}{2}$) indices larger than $n$},~~~~ \\
0&\mbox{otherwise}.\end{array}\right.$$
Then $\| \ST \| = 2^{m-1 \over 2}\| \Psi\|$.  Thus, we have
$$Q(\Psi) = Z(\ST) \le \rho(m, 2n)\| \ST\| = 2^{m-1 \over 2}\rho(m, 2n)\| \Psi \| \le {2^{m-1 \over 2}\rho(m, 2n) \over \mu(m, n)}Z(\Psi).$$
The conclusion follows.
\qed

As stated in the last section, this theorem initiates further study
on the QR-ratio.  what is its exact value?   Is it attainable?  
These remain for future research.

\section{Final remarks}

This paper reveals that the quantum eigenvalue problem has a close relationship with the Z-eigenvalue problem.    Two main results of this paper are Theorems \ref{rthm} and \ref{mthm}.   Theorem \ref{rthm} discusses the possibility to convert the problem for finding the entanglement eigenvalue of a real symmetric tensor to the Z-spectral radius problem of the same tensor, while Theorem \ref{mthm} converts the Q-eigenvalue problem of a complex symmetric tensor $\Psi$ to the Z-eigenvalue problem of a real symmetric tensor $\ST$, whose dimension is twice of the dimension of $\Psi$. The QR-ratio is introduced and proved to be finite.   Further research on these is needed.

{\bf Acknowledgment}  We are thankful to Shenglong Hu for his
comments.


\begin{thebibliography}{99}{\small}

\bibitem{ba} V. Balan, {\em Spectral properties and applications of numerical
multi-linear algebra of m-root structures}, in: Hypercomplex Numbers
in Geometry and Physics, ed., ``Mozet'', Russia, 2, 101-107, 2008.

\bibitem{ba09} V. Balan, {\em Numerical multilinear algebra of symmetric m-root
structures: Spectral properties and applications}, Symmetry: Culture
and Science, Symmetry Festival 2009, Symmetry in the History of
Science, Art and Technology; Part 2; Geometric Approaches to
Symmetry 2010, 21, 119-131, 2009.

\bibitem{ba11} V. Balan, {\em Spectra of symmetric tensors and m-root Finsler
models}, Linear Algebra and Its Applications, 436, 152-162, 2012.

\bibitem{ba10} V. Balan and N. Perminov, {\em Applications of resultants in the
spectral $M$-root framework}, Applied Sciences, 12, 20-29, 2010.

\bibitem{bv} L. Bloy and R. Verma, {\rm On computing the underlying fiber
directions from the diffusion orientation distribution function},
in: Medical Image Computing and Computer-Assisted Intervention --
MICCAI 2008,  D. Metaxas, L. Axel, G. Fichtinger and G. Sz\'ekeley,
eds., Springer-Verlag, Berlin, 1-8, 2008.

\bibitem{bh} D.C. Brody and L.P. Hughston, {\rm Geometric quantum mechanics}, J. Geom.
Phys., 38, 19--53, 2001.

\bibitem{cs} D. Cartwright and B. Sturmfels, {\em The number of eigenvalues of
a tensor}, to appear in: Linear Algebra and Its Applications.


\bibitem{cpz} K.C. Chang, K. Pearson and T. Zhang, {\em On eigenvalues of real
symmetric tensors}, J. Math. Anal. Appl., 350, 416-422, 2009.

\bibitem{cpz1} K.C. Chang, K. Pearson and T. Zhang, {\em Some variational principles of the $Z$-eigenvalues for nonnegative tensors}, School of
Mathematical Sciences, Peking University, December 2011.

\bibitem{cz} K.C. Chang and T. Zhang, {\em On the uniqueness and non-uniqueness of the $Z$-eigenvector for transition probability tensors}, School of
Mathematical Sciences, Peking University, March 2012

\bibitem{cxz} L. Chen, A. Xu and H. Zhu, {\em Computation of the geometric measure of
entanglement for pure multiqubit states}, Phys. Rev. A, 82, 032301,
2010.

\bibitem{clo} D. Cox, J. Little and D. O'Shea, {\sl Using Algebraic Geometry}, Springer-Verlag, New York, 1998.

\bibitem{gkz} I.M. Gelfand, M.M. Kapranov and A.V. Zelevinsky, {\sl
Discrimants, Resultants and Multidimensional Determinants},
Birkh\"auser, Boston, 1994.

\bibitem{hmmov} M. Hayashi, D.Markham, M. Murao, M. Owari and S. Virmani, {\rm The
geometric measure of entanglement for a symmetric pure state with
non-negative amplitudes}, J. Math. Phys., 50, 122104, 2009.

\bibitem{hs} J.J. Hilling and A. Sudbery, {\em The geometric measure of multipartite
entanglement and the singular values of a hypermatrix}, J. Math.
Phys., 51, 072102, 2010.

\bibitem{hq} S. Hu and L. Qi, {\em Algebraic connectivity of an even uniform
hypergraph}, to appear in: Journal of Combinatorial Optimization.

\bibitem{hqz} S. Hu, L. Qi and G. Zhang, The geometric measure of entanglement of
pure states with nonnegative amplitudes and the spectral theory of
nonnegative tensors, arXiv:1203.3675v5.

\bibitem{hkwg} R. H\"{u}bener, M. Kleinmann, T.-C. Wei, C.G. Guill\'{e}n, C.
Gonzalez-Guillen and O. G\"{u}hne, {\rm Geometric measure of
entanglement for symmetric states}, Phys. Rev. A, 80, 032324, 2009.

\bibitem{km} T.G. Kolda and J.R. Mayo, {\em Shifted power method for computing
tensor eigenpairs}, SIAM Journal on Matrix Analysis and
Applications, 32, 1095-1124, 2011.

\bibitem{lqz} A. Li, L. Qi and B. Zhang, {\em $E$-characteristic polynomials of tensors}, to appear in: Commun. Math. Sci.

\bibitem{lqy} G. Li, L. Qi and G. Yu, {\em The $Z$-eigenvalues of a symmetric
tensor and its application to spectral hypergraph theory},
Department of Applied Mathematics, University of New South Wales,
December 2011.

\bibitem{li} L-H. Lim, {\em Singular values and eigenvalues of tensors: A
variational approach}, Proceedings of the First IEEE International
Workshop on Computational Advances in Multi-Sensor Adaptive
Processing (CAMSAP), December 13-15, 129-132, 2005.

\bibitem{nqww} G. Ni, L. Qi, F. Wang and Y. Wang, {\em The degree of the
E-characteristic polynomial of an even order tensor}, Journal of
Mathematical Analysis and Applications, 329, 1218-1229, 2007.

\bibitem{nc} M.A. Nielsen and I.L. Chuang, {\em Quantum Computing and Quantum
Information}, Cambridge University Press, Cambridge, 2000.

\bibitem{or} J.M. Ortega and W.C. Rheinboldt, {\em Iterative Solution of Nonlinear Equations in Several Variables}, Academic Press, New York, 1970; Republication: SIAM, Philadelphia, 2000.

\bibitem{odv} R. Or\'{u}s, S. Dusuel and J. Vidal, {\em Equivalence of critical scaling
laws for many-body entanglement in the Lipkin-Meshkov-Glick Model},
Phys. Rev. Lett., 101, 025701, 2008.

\bibitem{qi} L. Qi, {\em Eigenvalues of a real supersymmetric tensor},
Journal of Symbolic Computation, 40, 1302-1324, 2005.

\bibitem{qi06} L. Qi, {\em Rank and eigenvalues of a supersymmetric tensor, a multivariate homogeneous
polynomial and an algebraic surface defined by them},  Journal of
Symbolic Computation, 41, 1309-1327, 2006.

\bibitem{qi07} L. Qi, {\em Eigenvalues and invariants of tensors}, Journal of
Mathematical Analysis and Applications, 325, 1363-1377, 2007.

\bibitem{qi11} L. Qi, {\em The best rank-one approximation ratio of a tensor space},
SIAM Journal on Matrix Analysis and Applications,  32, 430-442,
2011.

\bibitem{qi12}L. Qi, {\it The minimum Hartree value for the quantum entanglement problem,} arXiv: 2002.2983vl.

\bibitem{qww} L. Qi, F. Wang and Y. Wang, {\em Z-eigenvalue methods for a global
polynomial optimization problem}, Mathematical Programming, 118,
301-316, 2009.

\bibitem{qyw} L. Qi, G. Yu and E.X. Wu, {\em Higher order positive
semi-definite diffusion tensor imaging}, SIAM Journal on Imaging
Sciences, 3, 416-433, 2010.

\bibitem{rv} S. Ragnarsson and C.F. Van Loan, {\em Block tensors and
symmetric embeddings}, to appear in: Linear Algebra and Its
Applications.

\bibitem{wg} T.C. Wei and P.M. Goldbart, {\em Geometric measure of entanglement
and applications to bipartite and multipartite quantum states},
Phys. Rev. A, 68, 042307, 2003.

\bibitem{za} T. Zhang, {\em Existence of real eigenvalues of real tensor}, Nonlinear Analysis, 74, 2862-2868, 2011.

\bibitem{zqy} X. Zhang, L. Qi and Y. Ye, {\em The cubic spherical optimization problem},
{\em Mathematics of Computation} 81, 1513-1526, 2012.


\end{thebibliography}
\end{document}